\documentclass[
12pt,a4paper]{article}

\usepackage[marked,printedin,standardeqnum
]{myart2k}
\ppnum{LEEDS-MATH-PURE-2001-31}{\htmladdnormallink{math/9808040}%
{http://xxx.lanl.gov/abs/math/9808040}}{1998}

\input{mydef}
\titleshort{Polynomial Sequences of Binomial Type and Path Integrals}
\authorshort{Vladimir V. Kisil}

\newcommand{\column}[2]{\left(\hspace{-5pt}\begin{array}{c}#1\\#2
\end{array}\hspace{-5pt}\right)}

\begin{document} \title[Polynomial Sequences and Path
Integrals]{Polynomial Sequences of Binomial Type\\ and Path
  Integrals\thanks{Partially supported by the grant YSU081025 of
    \emph{Renessance} foundation (Ukraine).}}  \author[Vladimir V.
Kisil] {\href{http://maths.leeds.ac.uk/kisilv~/}{Vladimir V.
    Kisil}\footnote{On leave form the Odessa National University
    (Ukraine).}}
\address{School of Mathematics\\
  University of Leeds\\ Leeds LS2\,9JT\\ UK }
\email{\href{mailto:kisilv@maths.leeds.ac.uk}%
  {kisilv@maths.leeds.ac.uk}}
\urladdr{\href{http://maths.leeds.ac.uk/~kisilv/}%
  {http://maths.leeds.ac.uk/\~{}kisilv/}} 
\date{August 8, 1998; Revised October 19, 2001}
\maketitle
\begin{abstract}
  Polynomial sequences $p_n(x)$ of binomial type are a principal tool
  in the umbral calculus of enumerative combinatorics.  We express
  $p_n(x)$ as a \emph{path integral} in the ``phase space''
  $\Space{N}{} \times {[-\pi,\pi]}$.  The Hamiltonian is
  $h(\phi)=\sum_{n=0}^\infty p_n'(0)/n!  \, e^{in\phi}$ and it
  produces a Schr\"odinger type equation for $p_n(x)$.  This
  establishes a bridge between enumerative combinatorics and quantum
  field theory. It also provides an algorithm for parallel quantum
  computations.  
  \AMSMSC{05A40}{05A15, 58D30, 81Q30, 81R30, 81S40}
\end{abstract}
{\small\tableofcontents} \newpage \epigraph{Under the inspiring
  guidance of Feynman, a short-hand way of expressing---and of
  thinking about---these quantities have been developed.}{Lewis H.
  Ryder~\cite[Chap.~5]{Ryder96}.}{}
\section{Introduction on Polynomial Sequences of Binomial Type}
The umbral calculus~\cite{Foundationiii,RomRota78,Rota64a,KahOdlRota73}
in enumerative combinatorics uses polynomial sequences of binomial
type as a principal ingredient. A polynomial sequence $p_n(x)$ is said
to be of \emph{binomial type}~\cite[p.~102]{RomRota78} if $p_0(x)=1$
and if it satisfies the \emph{binomial identity},
\begin{eqnordisp}[eq:binomial]
  p_n(x+y)=\sum_{k=0}^n \column{n}{k} p_k(x) p_{n-k}(y)
\end{eqnordisp}
for all $n\in\Space{N}{}$ and $x$, $y\in\Space{R}{}$.
\begin{example}\label{ex:binomial}
  Among many examples~\cite[\S~5]{RomRota78} of polynomial sequences
  of binomial type we mention only very few:
  \begin{enumerate}
  \item The power monomials $p_n(x)=x^n$.
  \item The rising factorial sequence $p_n(x)=x(x+1) \cdots (x+n-1)$.
  \item The falling factorial sequence $p_n(x)=x(x-1) \cdots (x-n+1)$.
  \item The Abel polynomials $A_n(x)=x(x-an)^{n-1}$.
  \item The Laguerre polynomials
    \begin{eqnordisp}[eq:laguerre] L_n(x)=\sum_{k=0}^n \frac{n!}{k!}
      \column{n-1}{k-1} (-x)^k.
    \end{eqnordisp}
  \end{enumerate}
\end{example}

Using \emph{delta functionals} one can show
that~\cite[p.~106]{RomRota78} that for any polynomial sequence of
binomial type $p_n(x)$:
\begin{eqnordisp}[eq:a_n0]
  p_n(0)=\delta_{n0}=\left\{
    \begin{array}{ll}
      1, & n=0;\\
      0, & n>0.
    \end{array}
\right.
\end{eqnordisp}
For a polynomial sequence $p_n(x)$ we could consider its coefficients
$a_{n,k}$ in its decomposition over power monomials $x^k$:
\begin{eqnordisp}[]
  p_n(x) = \sum_{k=0}^n a_{n,k} x^k.
\end{eqnordisp}
The coefficients $a_{n,k}$ are an example of \emph{connecting
  constants}~\cite[\S~6]{RomRota78} between sequences $p_n(x)$ and
$x^n$. 
The definition~\eqref{eq:binomial} of polynomial sequences of binomial
type can be expressed in the term of coefficients
$a_{n,k}$ as follows~\cite[Prop.~4.3]{RomRota78}: 
\begin{eqnordisp}[eq:a-recurr]
  \column{i+j}{i} a_{n,i+j}= \sum_{k=0}^n \column{n}{k} a_{k,i}
  a_{n-k,j}.
\end{eqnordisp} Therefore these coefficients are highly
interdependent.  Particularly from~\eqref{eq:a_n0} it follows that
$a_{n,0}=\delta_{n,0}$. Moreover one can observe
from~\eqref{eq:a-recurr} that the knowledge of the sequence $c_n=a_{n,1}$,
$n=0,1,\ldots$ allows to reconstruct all $a_{n,k}$ by a recursion.
Sequence $c_n=a_{n,1}$ is known as \emph{cumulants} associated with
polynomials $p_n(x)$ and the related probability measure, see
e.g. \cite{Mattner99} for a survey.
One may state the following
\begin{prob}\label{pr:binom}
  Give a direct expression for a polynomial sequence $p_n(x)$ of
  binomial type from the sequence of its cumulants
  $c_n=a_{n,1}=p_n'(0)$, $n=0,1,\ldots$.
\end{prob}

There are several well known solutions to that problem. The umbral
calculus suggests the following
procedure\footnote{\label{fn:thanks-1}I am grateful to 
  the first anonymous referee for its indication.}
  (see~\cite{Knuth99a} for an
interesting exposition). Let us define the formal power series:
\begin{eqnordisp}[eq:cumulant-gen-funct]
  f(t)=\sum_{n=1}^\infty c_n \frac{t^n}{n!},
\end{eqnordisp}
then the function $e^{xf(t)}$ is a generating function for
polynomials $p_n(x)$:
\begin{eqnordisp}[eq:polynom-gen-func]
  \sum_{n=0}^\infty p_n(x) \frac{t^n}{n!}= e^{xf(t)}.
\end{eqnordisp}

The main goal of this paper is to present a different solution to the
Problem~\ref{pr:binom}. It is probably not the best from the
computation point of view unless computations are made in aprallel by
a quantum computer, see Remark~\ref{re:computing}. We will express
$p_n(x)$ as a \emph{path integral}~\eqref{eq:binom-path} in the
``phase space'' $\Space{N}{} \times [-\pi,\pi]$.  The Fourier
transform $h(\phi)=\sum_{n=0}^\infty \frac{c_{n}}{n!}e^{in\phi}$ of
the sequence $c_{n}/n!=p_n'(0)/n!$ plays the r\^ole of a
\emph{Hamiltonian}.  Thus we establish a connection between
enumerative combinatorics and quantum field theory.  This connection
could be useful in both areas.  It is interesting to note that quantum
mechanical commutation relations $PQ-QP=I$ was applied to umbral
calculus by \person{J.~Cigler}~\cite{Cigler78}.  Existence of a
``yet-to-be-discovered stochastic process'' connected with Laguerre
polynomials was conjectured by \person{S.M.~Roman} and
\person{G.-C.~Rota}~\cite{RomRota78}.  Note also that the approach to
the umbral calculus developed
in~\cite{Foundationiii,RomRota78,Rota64a,KahOdlRota73} can be a
prototype for a general construction of \emph{coherent states}
(\emph{wavelets}) in Banach spaces~\cite{Kisil98a}.

\section{Preliminaries on Path Integrals} 

Let a physical system has the \emph{configuration space}
$\Space{Q}{}$.  This means that we can label results of our
measurements performed on the system at any time $t_0 \in \Space{R}{}$
by points of $\Space{Q}{}$.  For a quantum system dynamics the
principal quantity is the \emph{propagator}
\cite[\S~2.2]{FeynHibbs65}, \cite[\S~5.1]{Ryder96}
$K(q_2,t_2;q_1,t_1)$---a complex valued function defined on
$\Space{Q}{} \times \Space{R}{} \times \Space{Q}{} \times
\Space{R}{}$.  It is a probability amplitude for a transition $q_1
\rightarrow q_2$ from a state $q_1$ at time $t_1$ to $q_2$ at time
$t_2$.  The \emph{density of probability} of this transition is
postulated to be
\begin{eqnordisp}[]
  P(q_2,t_2;q_1,t_1) = \modulus{K(q_2,t_2;q_1,t_1)}^2.
\end{eqnordisp} 

The fundamental assumption about the quantum world is the
\emph{absence of trajectories} for a system's evolution through the
configurational space $\Space{Q}{}$: the system at any time $t_i$
could be found at any point $q_i$.  Then the transition amplitude $q_1
\rightarrow q_2$ is a result of all possible transitions amplitudes
$q_1 \rightarrow q_i \rightarrow q_2$ integrated over $\Space{Q}{}$.

Mathematically this can expressed in the term of propagator
$K$ as follows. Let us fix any $t_i$,
$t_1< t_i < t_2$.  There exists a measure $dq$ on $\Space{Q}{}$ such that
\begin{eqnordisp}[eq:k-proj]
  K(q_2,t_2;q_1,t_1)= \int_{\Space{Q}{}} K(q_2,t_2;q_i,t_i)
  K(q_i,t_i;q_1,t_1)\, dq_i.
\end{eqnordisp}
This identity could be thought as a \emph{recurrence formula} for the
propagator $K$.

Another assumption on propagator is its skew
symmetry
\begin{eqnordisp}[eq:k-adjoint]
  K(q_2,t_2;q_1,t_1)= \overline{K(q_1,t_1;q_2,t_2)},
\end{eqnordisp} it ensures that all transitions $q_1 \rightarrow q_2$
are done by unitary operators.

It is useful to consider $K(q_2,t_2;q_1,t_1)$ as a kernel for the
associated integral operator $K: \FSpace{L}{2}(\Space{Q}{},dq)
\rightarrow \FSpace{L}{2}(\Space{Q}{},dq)$ defined as
\begin{eqnordisp}[eq:k-int]
  [K \phi](q,t) = \int_{\Space{Q}{}} K(q,t;q_i,t_i)
  \phi(q_i,t_i)\,dq_i,
\end{eqnordisp}
with $t_i$ is fixed.  Operator $K$ is defined on functions
$\phi(q,t)\in\FSpace{L}{2}(\Space{Q}{},dq)$ such that the right-hand
side of~\eqref{eq:k-int} does not depend on $t_i$.  Then
\eqref{eq:k-proj} is read as $K^2=K$ and \eqref{eq:k-adjoint} implies
$K^*=K$.  By other word $K$ is an orthoprojection\footnote{For a model of
  relativistic quantisation generated by two orthoprojections
  see~\cite{Kisil95b}.} from $\FSpace{L}{2}(\Space{Q}{},dq)$ to its
proper subspace $\object{Im}K$.

Another consequence of~\eqref{eq:k-proj}:  functions
$\phi(g_2,t_2)$ which belongs to image of operator $K$ are fixed
points for $K$:
\begin{eqnordisp}[eq:wave]
  \phi(q,t)=[K \phi](q,t) = \int_{\Space{Q}{}} K(q,t;q_i,t_i)
  \phi(q_i,t_i)\,dq_i.
\end{eqnordisp} These are \emph{wave functions} and describe actual
\emph{states} of the system.  It follows from~\eqref{eq:k-proj} that
all functions $\kappa(q,t)=K(q,t;q_0,t_0)$ for arbitrary fixed $q_0\in
\Space{Q}{}$, $t_0\in \Space{R}{}$ are wave functions too.  The function
$\modulus{\phi(q,t)}^2$ gives the probability distribution over
$\Space{Q}{}$ to observe the system in a point $q$ at time $t$ if it
was in $q_0$ at time $t_0$.

\person{R.~Feynman} developing ideas of \person{A.~Einstein},
\person{Smoluhovski} and \person{P.A.M.~Dirac} proposed an expression
for the propagator via the ``integral over all possible paths''. In a
different form it can be written as follows~\cite[(5.13)]{Ryder96}
\begin{eqnordisp}[eq:path]
  K(q_2,t_2;q_1,t_1)=\int \frac{\mathcal{D}q\, \mathcal{D}p}{h} \exp
  \left( \frac{i}{\hbar} \int\limits_{t_1}^{t_2} dt \left(p\dot{q}-H(p,q)
    \right) \right).
\end{eqnordisp} Here $H(p,q)$ is the classical
Hamiltonian~\cite[Chap.~9]{Arnold91} of the system---a function
defined over the \emph{phase space}. In the simplest case the phase
space is the Cartesian product of the configuration space
$\Space{Q}{}$ and the space of momenta dual to it.  Points in phase
space have a double set of coordinates $(q,p)$, where $q$ is points of
the configuration space $\Space{Q}{}$ and $p$ is a coordinate in the
space of momenta.  The inner integral in~\eqref{eq:path} is taken over
a path in the phase space. The outer integral is taken over ``all
possible paths with respect to a measure $\mathcal{D}q\, \mathcal{D}p$
on paths in the phase space''. We construct an example of path
integral~\eqref{eq:N-approx}--\eqref{eq:path-form}) similar
to~\eqref{eq:path} in the next section.

From reading of physics textbooks~\cite{FeynHibbs65,Ryder96} one can
obtain an impression that the measure $\mathcal{D}q\, \mathcal{D}p$
could be well defined from an unspecified complicated mathematical
construction\footnote{The editor preface to the Russian translation
  of~\cite{FeynHibbs65} tells: ``There is no a rigorous definition of
  the path integral in the book of Feynman and Hibbs\ldots However it
  is not very important for a physicist in most cases; he only need an
  assurance that a rigorous proof could be obtained''.  The recent
  textbook~\cite{Ryder96} which were already reprinted 9 times in two
  editions refereed for a rigorous path integral construction to the
  old paper~\cite{GelfandYaglom60a}, however it is known that this
  approach is not complete, see e.g. the review of this article
  \MR{17,1261c} as well as \cite{Daletski62a} and
  \cite[Notes~X.11]{SimonReed75}.}.
But this is \emph{not} true up to now: \emph{there is no a general
mathematically rigorous definition of the measure in the path space}.
However there are different approaches of regularisation of path
integrals, see~\cite{Dynin98,ShabanovKlauder98} for recent accounts.

It is known that path integrals give fundamental solutions of
evolutionary equations~\cite{Daletski62a}. Particularly wave
functions~\eqref{eq:wave} are solutions to the Schr\"odinger equation
\cite[\S~5.2]{Ryder96}, \cite[\S~4.1]{FeynHibbs65}
\begin{eqnordisp}[eq:schrodinger]
  i\hbar \frac{\partial \phi(q,t)}{\partial t}= \widehat{H}(q,p)
  \phi(q,t),
\end{eqnordisp} where $\widehat{H}(q,p)$ is an operator which is the
quantum Hamiltonian.  It was mentioned that functions
$\kappa(q,t)=K(q,t;q_0,t_0)$ are wave functions and thus are solutions
to~\eqref{eq:schrodinger}.  In fact it follows from~\eqref{eq:wave}
that they are \emph{fundamental solutions}, i.e. any other solution is
a linear combination (\emph{superposition} in physical language) of
those.  Therefore the propagator $K(q_2,t_2;q_i,t_i)$ contains the
complete information on dynamics of a quantum system defined by the
equation~\eqref{eq:schrodinger}.

\section{Polynomial Sequences from Path Integrals} 

Let the configuration space $\Space{Q}{}$ be a semigroup with an
operation denoted by $+$ and let a propagator $K(q_2,t_2;q_1,t_1)$ be
homogeneous in time and space, i.e.
\begin{eqnordisp}[eq:homogeneuity]
  K(q_2+q,t_2+t;q_1+q,t_1+t)=K(q_2,t_2;q_1,t_1),
\end{eqnordisp}
for all $q, q_1, q_2\in \Space{Q}{}$, $t, t_1, t_2 \in \Space{R}{}$.
Such a propagator corresponds to a free time-shift invariant system.
Then $K(q_2,t_2;q_1,t_1)$ is completely defined by values of the
function $S(q,t)$ such that
$S(q_2-q_1,t_2-t_1)=K(q_2,t_2;q_1,t_1)$. From~\eqref{eq:k-proj} we
can deduce a ``recurrence'' identity for the function $S(q,t)$
\cite{Daletski62a}:
\begin{eqnordisp}[eq:token]
  S(q,t_1+t_2)= \int_{\Space{Q}{}} S(q_i,t_1) S(q-q_i,t_2)\,dq_i.
\end{eqnordisp}

An equation of the form~\eqref{eq:token} was taken in~\cite{Kisil97b}
as the definition of \emph{tokens} between two cancellative semigroups
$\Space{Q}{}$ and $\Space{R}{+}$. It allows to express a time shift by
the integral of all space shifts over $\Space{Q}{}$. Tokens are
closest relatives to \emph{sectional coefficients} considered by
Henle~\cite{Henle75}.  A polynomial sequence $p_n(x)$ of binomial
type is another example of token~\cite{Kisil97b}: the sequence
$q_n(x)=p_n(x)/n!$ is a token between cancellative semigroup
$\Space{N}{}$ and $\Space{R}{}$, i.e. from~\eqref{eq:binomial} follows
a realization of~\eqref{eq:token} in the form
\begin{eqnordisp}[eq:q-polyn]
  q_n(x+y)=\sum_{k=0}^\infty q_{k}(y) q_{n-k}(x),
\end{eqnordisp}
i.e. in this case the configuration space $\Space{Q}{}$ is the
discrete set of nonnegative numbers and the integration reduces to
summation. 

If polynomial sequence $q_n(x)$ represents a particular propagator for
a configuration space $\Space{Q}{}=\Space{N}{}$ then one may wish to
obtain a path integral representation in the form~\eqref{eq:path} for
it.  We do it in a way parallel to the path integral deduction in
\cite[\S~5.1]{Ryder96}.

Let $x\in \Space{R}{}$ be fixed and $N\in \Space{N}{}$ be an arbitrary
positive integer. Then repetitive use of~\eqref{eq:q-polyn} allows to
write
\begin{eqnordisp}[eq:sum-iter]
  q_n(x)= \relstack{ \sum \sum \cdots \sum}{ k_1+k_2+\cdots+k_N=n}
  q_{k_1}\!\left(\frac{x}{N} \right)\, q_{k_2}\! \left(\frac{x}{N}
  \right)\cdots q_{k_N}\! \left(\frac{x}{N} \right).
\end{eqnordisp} The first two terms in the Taylor expansion of
$q_{k_j}$, $j=1,\ldots,N$ are:
\begin{eqnordisp}[eq:taylor-exp]
  q_{k_j}\left(\frac{x}{N}\right) = \delta_{k_j 0} + q_{k_j}'(0)
  \frac{x}{N} + o\! \left(\frac{x}{N}\right),
\end{eqnordisp}
where the Kronecker delta as the first term come from~\eqref{eq:a_n0}.
We use the Pontrjagin duality between $\Space{Z}{}$ and the unit disk
$\Space{T}{}=[-\pi,\pi]$ to construct an integral
resolution of the Kronecker delta:
\begin{eqnordisp}[eq:delta]
  \delta_{k 0}= \frac{1}{2\pi} \int\limits_{-\pi}^\pi e^{-ip k}\,dp.
\end{eqnordisp} It is also convenient to introduce the Fourier
transform $h(p)$ for the numerical sequence $q_{k}'(0)$ from the same
Pontrjagin duality:
\begin{eqnordisp}[eq:hamiltonian]
  h(p)=\sum_{k=0}^\infty q_k'(0) e^{ipk}
  \qquad\iff\qquad
  q_k'(0)= \frac{1}{2\pi} \int\limits_{-\pi}^\pi e^{-ipk}h(p)\,dp.
\end{eqnordisp}
Then we can transform the Taylor expansion~\eqref{eq:taylor-exp} as
follows: 
\begin{eqnarray}
q_{k}\left(\frac{x}{N}\right) & = & \delta_{k 0} + q_{k}'(0) \frac{x}{N} 
+ o\!\left(\frac{x}{N}\right) \nonumber\\
& = & \frac{1}{2\pi}\int\limits_{-\pi}^\pi e^{- ipk}dp + \frac{1}{2\pi}
\int\limits_{-\pi}^\pi e^{-ipk} h(p)\,dp \frac{x}{N}
+ o\!\left(\frac{x}{N}\right) \nonumber\\
& = & \frac{1}{2\pi}\int\limits_{-\pi}^\pi e^{- i p k}\left( 1+
 h(p) \frac{x}{N} \right)dp + o\!\left(\frac{x}{N}\right) \nonumber\\
&=& \frac{1}{2\pi} \int\limits_{-\pi}^\pi \exp\left(-i p k + h(p) \frac{x}{N}
\right) dp + o\!\left(\frac{x}{N}\right), \label{eq:q-delta}
\end{eqnarray}
where the last transformation is based on the estimation 
\begin{eqnordisp}
  \exp\left(h(p) \frac{x}{N}\right)=  \left(1+
 h(p) \frac{x}{N}\right) + o\!\left(\frac{x}{N}\right).
\end{eqnordisp}

Substituting $N$ copies of~\eqref{eq:q-delta} to~\eqref{eq:sum-iter}
we obtain
\begin{eqnarray}
\lefteqn{q_n(x)  =  \relstack{ \sum \sum \cdots \sum}{ k_1+k_2+\cdots+k_N=n} 
q_{k_1}\!\left(\frac{x}{N} \right)\, q_{k_2}\! \left(\frac{x}{N} 
\right)\cdots q_{k_N}\! \left(\frac{x}{N} \right) \nonumber}\\
& = & \relstack{ \sum \cdots \sum}{ k_1+\cdots+k_N=n}
\left( \frac{1}{(2\pi)^N} \int\limits_{-\pi}^\pi \exp\left(-i p_1 k_1 + h(p_1)
\frac{x}{N}\right) dp_1 \right.\nonumber\\
&& \qquad\qquad \left. \times \cdots \times
\int\limits_{-\pi}^\pi \exp\left(-i p_N k_N + h(p_N) \frac{x}{N}\right) dp_N 
+ o\!\left(\frac{x}{N}\right) \right) \nonumber\\
& = & \relstack{ \sum \cdots \sum}{ k_1+\cdots+k_N=n} \left(
\frac{1}{(2\pi)^N} \int\limits_{-\pi}^\pi \cdots
\int\limits_{-\pi}^\pi 
\exp\left(-i p_1 k_1 + h(p_1) \frac{x}{N}\right) \right.\nonumber\\ && \qquad
\qquad \left.\times\cdots \times
 \exp\left(-i p_N k_N + h(p_N) \frac{x}{N}\right) dp_1 \ldots dp_N
+ o\!\left(\frac{x}{N}\right) \right) \nonumber\\
& = & \relstack{ \sum \cdots \sum}{ k_1+\cdots+k_N=n} \left(
\frac{1}{(2\pi)^N} \int\limits_{-\pi}^\pi \cdots
\int\limits_{-\pi}^\pi 
\exp\left({\sum_{l=1}^n \left(-i p_l k_l + h(p_l) \frac{x}{N}\right)}\right)
 dp_1 \ldots dp_N \right.  \nonumber\\
&& \qquad \qquad \left.+ o\!\left(\frac{x}{N}\right) \right)
\label{eq:N-approx}
\end{eqnarray}
Now we can observe that the expression $\sum_{l=1}^n \left(-i p_l
  k_l + h(p_l) \frac{x}{N}\right)$ looks like an integral sum for the
function $-ip(t)k'(t)+h(p(t))$ on the interval $t\in[0,x]$.

At this point we say together with physicists~\cite[\S~5.1]{Ryder96}
magic words and \eqref{eq:N-approx} ``could be symbolically written
for $N\rightarrow \infty$'' as
\begin{eqnordisp}[eq:path-form]
  q_n(x)= \int\!\mathcal{D}k\mathcal{D}p\, \exp\!\!\int\limits_0^x
    \left(-ip(t)k'(t)+h(p(t))\right)\, dt,
\end{eqnordisp} where the first integration is taken over all possible
paths $k(t):[0,x] \rightarrow \frac{1}{x}\Space{N}{}$, such that
$k(0)=0$ and $k(x)=\frac{n}{x}$, and the path $p(t): [0,x]\rightarrow
[-{\pi},{\pi}]$ is unrestricted. It is enough to consider only paths
$k(t)$ with monotonic grow---other paths make the zero contributions.
Indeed the identity $q_l(x)\equiv 0$ for $l<0$ (made by an agreement)
implies that contribution of all paths with $k'(t)<0$ at some point
$t$ vanish. Here (and in the inner integral of~\eqref{eq:path-form})
$k'(t)$ means the derivative of the path $k(t)$ in the distributional
sense, i.e. it is the Dirac delta function times $\frac{j}{x}$ in the
points where $k(t)$ jumps from one integer (mod $\frac{1}{x}$) value
$k(t-0)=\frac{m}{x}$ to another $k(t+0)=\frac{m+j}{x}$.

In fact there is no any magic in the above path integral. Due to the
property $q_n(x) \equiv 0$ for $n<0$ paths $k(t)$ take values on the
finite set $k_0$, $k_1$,\ldots, $k_n$, otherwise they make the zero
contribution to the path integral as was mentioned above. Therefore
the set of paths $k(t)$ is countable. Under such a condition the
situation is equivalent to a finite Markov process with a well-defined
measure on the path space according to Kolmogorov theorem, see for
example~\cite[\S~1]{Daletski62a} or \cite[Th.~2.1]{Simon79}. Moreover
at the beginning of the next section we show how to calculate the
path integral~\eqref{eq:path-form} as an ordinary integral.

Now we can give a solution to Problem~\ref{pr:binom} using the
relation $p_n(x)=n!q_n(x)$ between $p_n(x)$ and $q_n(x)$:
\begin{eqnordisp}[eq:binom-path]
  p_n(x)= n!\int\!\mathcal{D}k\mathcal{D}p\,
  \exp\!\!\int\limits_0^x
  (-ipk'+h(p))\,dt, \qquad h(p)=\sum_{k=0}^\infty p_k'(0)
  \frac{e^{ipk}}{k!}. 
\end{eqnordisp} Note that a comparison of our Hamiltonian $h(p)$ with
the formal power series $f(t)$~\eqref{eq:cumulant-gen-funct} made out
of cumulants shows that
$
  h(p)=f(e^{ip}).
$
\begin{example} \label{ex:hamiltonian}
  Here is a list of Hamiltonians producing the polynomial sequences
  from Example~\ref{ex:binomial}.
\begin{enumerate}
\item The power monomials $p_n(x)=x^n$ are produced by $h(p)=e^{ip}$.
\item The rising factorial sequence $p_n(x)=x(x+1) \cdots (x+n-1)$ is
  produced by the Heaviside function
  $
    h(p)=\chi(p)=\sum_{k=1}^\infty \frac{e^{ipk}}{k}.
  $
\item The falling factorial sequence $p_n(x)=x(x-1) \cdots (x-n+1)$ is
  produced by the shifted Heaviside function
  \begin{eqnordisp}
    h(p)=\chi(p-\pi)=\sum_{k=1}^\infty \frac{(-1)^k e^{ipk}}{k}
    =\sum_{k=1}^\infty \frac{e^{i(p-\pi)k}}{k} .
  \end{eqnordisp}
\item The Abel polynomials $A_n(x)=x(x-an)^{n-1}$ are produced by
  $h(p)=\sum_{k=1}^\infty \frac{(ak)^{k-1}}{k!} e^{ipk}$.
\item The Laguerre polynomials \eqref{eq:laguerre} are generated by
  the delta function 
  $
    h(p)=\delta(p)=\sum_{k=1}^\infty e^{ipk}.
  $
\end{enumerate}
\end{example}

\section{Some Applications}
\label{sec:schr-type-equat}

Let us start from the question\footnote{\label{fn:second-thanks}I am
grateful to the second anonymous referee who suggested this question.}
of calculational applicability of the path integral
formula~\eqref{eq:binom-path}. We transform this path integral
into an ordinary integration. Because we could not make any seriously
restricting assumptions about the function $h(p)$ the following
arguments have an algebraical-combinatorial nature rather than
analytic. They are similar to heuristic manipulations with formal
power series which are used in a deduction of
formula~~\eqref{eq:polynom-gen-func}.

For a convenience we scale by $\frac{1}{x}$ the parameter
$t$ in the inner integral of~\eqref{eq:path-form}:
\begin{eqnordisp}
  \int\limits_0^x \left(-ip(t)k'(t)+h(p(t))\right)\, dt=
  \int\limits_0^1 \left(-ip_1(t_1)k_1'(t_1)+xh(p_1(t_1))\right)\, dt_1,
\end{eqnordisp}
where on the right-hand side new
paths are $k_1(t_1):[0,1] \rightarrow \Space{N}{}$ with the endpoint
$k_1(0)=0$ and $k_1(x)={n}$, and the path $p_1(t_1): [0,1]\rightarrow
[-{\pi},{\pi}]$ is unrestricted. We continue renaming variables of
integration $t_1$, $k_1$, and $p_1$ by $t$, $k$ and $p$ respectively. 

As was mentioned before $k(t)$ is an increasing step functions
uniquely defined by a collection of $n$ points $t_j\in[0,1]$ (some of
them could coincide) where it jumps by $1$.  Thus $k'(t)=\sum_1^n
\delta (t-t_j)$ and consequently:
\begin{eqnordisp}
  \int\limits_0^1  -ip(t)k'(t)\,dt
  = -i\sum_{j=1}^n p(t_j).
\end{eqnordisp}
Thereafter we could transform~\eqref{eq:binom-path} as follows:
\begin{eqnarray}
p_n (x)& = & n! \int\!\mathcal{D}k\mathcal{D}p\,
  \exp\!\!\int\limits_0^1(-ip(t)k'(t)+xh(p(t)))\,dt \nonumber\\
& = & n! \int\!\mathcal{D}k\mathcal{D}p\,
  \exp \left(-i\sum_{j=1}^n
    p(t_j)+\int\limits_0^1 xh(p(t))\,dt\right) 
  \label{eq:path-transform-1}
\end{eqnarray}

Now we use the rotational symmetry of the unit circle $[-\pi,\pi]$ and
that paths $p(t)$ are completely unrestricted, therefore the whole set
$P$ of such paths is also rotational invariant. Thus the integration
over that set $P$ acts just like an averaging over the unit circle
$[-\pi,\pi]$. On the other hand a random path $k(t)$ is defined by a
random collection of $n$ points $t_j\in[0,1]$ consequently the
integration of $p(t_j)$ for any $j$ over all paths $k(t)$ will produce
$n$ independent equal uniformly distributed random variables on the
unit circle.  Thus we could finally express our path
integral~\eqref{eq:path-transform-1} as a simple integration:
\begin{eqnarray}
p_n (x) & = & n! \int\!\mathcal{D}k\mathcal{D}p\,
  \exp \left(-i\sum_{j=1}^n
    p(t_j)+\int\limits_0^1xh(p(t))\,dt\right) \nonumber\\
  & = & n!\int\limits_{-\pi}^{\pi}\exp
  \left(-i\sum_{j=1}^n    p+xh(p)\right)\,dp \nonumber\\
  & = & n!\int\limits_{-\pi}^{\pi}e^{-inp} e^{xh(p)}\,dp. 
  \label{eq:generating-func-2}
\end{eqnarray}

If we recall the connection $ h(p)=f(e^{ip})$ between our Hamiltonian
$h(p)$ and the generating function of cumulants
$f(t)$~\eqref{eq:cumulant-gen-funct} then the last
formula~\eqref{eq:generating-func-2} obtains the very simple meaning:
\emph{it uses the Fourier transform on $[-\pi,\pi]$ to extract $n$-th
  term out of the generating function~\eqref{eq:polynom-gen-func} for
  the polynomial sequence $p_n(x)$.}

\begin{rem}
  \label{re:computing} As we see that our
  formula~\eqref{eq:binom-path} can offer in conventional calculations
  the same (or slightly longer) procedure as the well-known umbral
  approach~\eqref{eq:polynom-gen-func}. The situation could change if
  hypothetical \emph{quantum computers}~\cite{Manin00a} will be able
  to calculate quantum propagators \emph{directly in parallel
    processes}. In this case the path integral
  formula~\eqref{eq:binom-path} should have advantages over the
  generating function~\eqref{eq:polynom-gen-func} where calculations
  \emph{must} be done in a row.
\end{rem}

Finally we consider a realisation of the Schr\"odinger
equation~\eqref{eq:schrodinger}. 
Homogeneous propagator $S(q,t)$~\eqref{eq:token} are wave
functions~\eqref{eq:wave} itself. Thus they should satisfy to the
Schr\"odinger like equation~\eqref{eq:schrodinger}. For polynomial
sequences $q_n(x)$ this equation takes the form
\begin{eqnordisp}[eq:schrod-form]
  \frac{\partial}{\partial x} q_n(x) = \widehat{H} q_n(x),
\end{eqnordisp}
where the operator $\widehat{H}$ is of the pseudodifferential
type~\cite{Shubin87,MTaylor81}
\begin{eqnarray}
\widehat{H} a_n & = & \fourier{p \rightarrow n}\, h(p) \fourier{n 
\rightarrow p}\, a_n \nonumber\\
& = & \int\limits_{-\pi}^\pi e^{-inp} h(p) \sum_{k=0}^\infty a_n
e^{inp} dp, \nonumber\\ &=& \sum_{k=0}^\infty h_k a_{n-k},
\label{eq:schrod-conv}
\end{eqnarray}
where $\fourier{p \rightarrow n}$ and $\fourier{n 
\rightarrow p}$ the Fourier transform in the indicated variables. 
The function $h(p)$ above is the Hamiltonian from the path
integral~\eqref{eq:path-form} and
\begin{eqnordisp}[]
  h_n= \frac{1}{2\pi}\int\limits_{-\pi}^\pi e^{-ipn} h(p)\,dp.
\end{eqnordisp} It turns to be just a convolution on $\Space{N}{}$
because Hamiltonian $h(p)$ depends only from $p$ and is independent
from $k$. The equation~\eqref{eq:schrod-conv} express the
property~\cite{Kisil97b} of tokens to intertwine shift invariant
operators between two cancellative semigroups. In the present case the
operators are the derivative with respect $x$ (i.e.  the convolution
with $\delta'(x)$ on $\Space{R}{}$) and the convolution with $h_n$ on
$\Space{N}{}$.

Let $\widehat{M}$ be the operator on sequences $\widehat{M}: \{a_n\}
\mapsto\{ a_n / n!\} $. Then $q_n(x)= \widehat{M} p_n(x) $ and because
$\widehat{M}$ and $\frac{\partial}{\partial x}$ commute we obtain for
a polynomial sequence $p_n(x)$ of binomial type the equation
\begin{eqnordisp}[eq:p-schrod]
  \frac{\partial}{\partial x} p_n(x)= \widehat{M}^{-1} \widehat{H}
  \widehat{M} p_n(x) \quad \iff \quad \frac{\partial}{\partial x}
  p_n(x)= n!  \sum_{k=0}^\infty h_k \frac{p_{n-k}(x)}{(n-k)!}.
\end{eqnordisp}
Of course the above equation follows from the differentiation of
identity~\eqref{eq:binomial} with respect to $y$ at the point
$y=0$. While this formula and the next example are rather elementary
we emphasise their relations to the quantum mechanical framework.  
\begin{example} \label{ex:schrodinger}
  The polynomial sequences from Examples~\ref{ex:binomial}
  and~\ref{ex:hamiltonian} satisfy to the following realizations of
  equation~\eqref{eq:p-schrod}:
\begin{enumerate}
\item The power monomials $p_n(x)=x^n$ satisfy to
\begin{eqnordisp}[] p_n'(x)= n! \sum_{k=0}^\infty \delta_{k1}
  \frac{p_{n-k}(x)}{(n-k)!}= \frac{n!}{(n-1)!} p_{n-1}(x)= n
  p_{n-1}(x).
\end{eqnordisp}
\item The rising factorial sequence $p_n(x)=x(x+1) \cdots (x+n-1)$
  satisfies to
\begin{eqnordisp}[]
  p_n'(x)= n! \sum_{k=0}^\infty \frac{1}{k} \frac{p_{n-k}(x)}{(n-k)!}.
\end{eqnordisp}
\item The falling factorial sequence $p_n(x)=x(x-1) \cdots (x-n+1)$
  satisfies to
\begin{eqnordisp}[]
  p_n'(x)= n! \sum_{k=0}^\infty \frac{(-1)^k}{k}
  \frac{p_{n-k}(x)}{(n-k)!}.
\end{eqnordisp}
\item The Abel polynomials $A_n(x)=x(x-an)^{n-1}$ satisfy to
\begin{eqnordisp}[] p_n'(x)= n! \sum_{k=0}^\infty
  \frac{(ak)^{k-1}}{k!}  \frac{p_{n-k}(x)}{(n-k)!} = \sum_{k=0}^\infty
  \column{n}{k} (ak)^{k-1} p_{n-k}(x) .
\end{eqnordisp}
\item The Laguerre polynomials \eqref{eq:laguerre} satisfy to
\begin{eqnordisp}[] p_n'(x)= n! \sum_{k=0}^\infty
  \frac{p_{n-k}(x)}{(n-k)!}.
\end{eqnordisp}
\end{enumerate}
\end{example}
\begin{rem}
  In our consideration we used primary the definition of
  token~\eqref{eq:token} and the Pontrjagin duality
  in~\eqref{eq:delta}--\eqref{eq:hamiltonian}.  Thus
  formula~\eqref{eq:path-form} has the same form for many other
  tokens~\cite{Kisil01b}, which we shall not consider here however.
  It follows from general properties of path integrals that for any
  admissible function $h(p)$ formula \eqref{eq:path-form} gives a
  token $q_n(x)$ and \eqref{eq:binom-path} produce a polynomial
  sequence $p_n(x)$ of binomial type.

  We obtain further generalisation if drop the
  homogeneity assumption~\eqref{eq:homogeneuity}. Let combinatorial
  quantities $p_{n,k}(x,y)$ where $n$, $k\in\Space{N}{}$ and $x$,
  $y\in\Space{R}{}$ be polynomials in the variable
  $x$ of the degree $n$ and in the variable
  $y$ of the degree $k$. We assume that they satisfy the identities
  \begin{eqnordisp}
    p_{n,k}(x,y)= \sum_{m=0}^\infty p_{n,m}(x,y)\,p_{m,k}(y,z),
  \end{eqnordisp} for any fixed $y$ such that $x\leq y\leq
  z$. Obviously this is an inhomogeneous version of the
  identities~\eqref{eq:binomial} and~\eqref{eq:token}. The above
  quantum mechanical derivation of path integral~\eqref{eq:path-form}
  can be carried (under some necessary assumptions) also in this
  case. The main difference in the resulting formulae is that the
  Hamiltonian $h(p,k)$ will depend from both variables
  $p\in[-\pi,\pi]$ and $k\in\Space{N}{}$. Consequently the
  Scr\"odinger equation~\eqref{eq:p-schrod} will not be any more a
  trivial convolution over $\Space{N}{}$. We left consideration of
  this case for some further papers.
\end{rem}

\emph{In the conclusion}: the main result of this paper is yet another
illustration to the old observation that mathematics is indivisible.

\section*{Acknowledgements} 

This work was partially supported by the grant YSU081025 of
\emph{Renessance} foundation (Ukraine). I am in debt to late Professor
Gian-Carlo Rota for introducing me to the world of umbral calculus.  I
am also grateful to two anonymous referees for their constructive
critique of the first version of this paper and many valuable
suggestions (particularly as indicated in the
footnotes~\ref{fn:thanks-1} and~\ref{fn:second-thanks}).

\small \bibliographystyle{plain}
\bibliography{abbrevmr,akisil,analyse,aphysics,acombin}

\newcommand{\noopsort}[1]{} \newcommand{\printfirst}[2]{#1}
  \newcommand{\singleletter}[1]{#1} \newcommand{\switchargs}[2]{#2#1}
  \newcommand{\irm}{\textup{I}} \newcommand{\iirm}{\textup{II}}
  \newcommand{\vrm}{\textup{V}} \providecommand{\cprime}{'}
\begin{thebibliography}{10}

\bibitem{Arnold91}
Vladimir~I. Arnold.
\newblock {\em Mathematical Methods of Classic Mechanics}.
\newblock Springer-Verlag, Berlin, {\noopsort{}}1991.

\bibitem{Cigler78}
J.~Cigler.
\newblock Some remarks on {Rota's} umbral calculus.
\newblock {\em Nederl. Akad. Wetensch. Proc. Ser. A}, 81:27--42, 1978.

\bibitem{Daletski62a}
Yuri~A. Daletski.
\newblock Path integrals connected with operator evolutionary equations.
\newblock {\em Uspehi Mat. Nauk}, 17(5(107)):3--115, 1962.

\bibitem{Dynin98}
Alexander Dynin.
\newblock A rigorous path integral construction in any dimension.
\newblock {\em Lett. Math. Phys.}, 44:317--329, 1998.
\newblock \eprint{http://xxx.lanl.gov/abs/math/9802058/}{math/9802058}.

\bibitem{FeynHibbs65}
R.P. Feynman and A.R. Hibbs.
\newblock {\em Quantum Mechanics and Path Integral}.
\newblock McGraw-Hill Book Company, New York, {\noopsort{}}1965.

\bibitem{GelfandYaglom60a}
I.~M. Gel{\cprime}fand and A.~M. Jaglom.
\newblock Integration in functional spaces and its applications in quantum
  physics.
\newblock {\em J. Mathematical Phys.}, 1:48--69, 1960.
\newblock \MR{17,1261c}.

\bibitem{Henle75}
Michael Henle.
\newblock
  \href{http://www.ams.org/leavingmsn?url=http://links.jstor.org/sici?sici=000%
2\%2D9947\%28197502\%29202\%3C1\%3ABEOD\%3E2.0.CO\%3B2\%2DQ\%26origin=MSN}{Bin%
omial enumeration on dissects}.
\newblock {\em Trans. Amer. Math. Soc.}, 202:1--39, 1975.
\newblock \MR{MR50:9601}.

\bibitem{Kisil01b}
Vladimir~V. Kisil.
\newblock Tokens: An algebraic construction common in combinatorics, analysis,
  and physics.
\newblock (In preparation).

\bibitem{Kisil95b}
Vladimir~V. Kisil.
\newblock Relativistic quantization and improved equation for a free
  relativistic particle.
\newblock {\em Phys. Essays}, 11(1):69--80, 1998.
\newblock \eprint{http://arXiv.org/abs/quant-ph/9502022/}{quant-ph/9502022}.
  \MR{99c:81117}.

\bibitem{Kisil98a}
Vladimir~V. Kisil.
\newblock Wavelets in {Banach} spaces.
\newblock {\em Acta Appl. Math.}, 59(1):79--109, 1999.
\newblock \eprint{http://arXiv.org/abs/math/9807141/}{math/9807141},
  \MR{2001c:43013}.

\bibitem{Kisil97b}
Vladimir~V. Kisil.
\newblock The umbral calculus: a model from convoloids.
\newblock {\em Z. Anal. Anwendungen}, 19(2):315--338, 2000.
\newblock \eprint{http://arXiv.org/abs/funct-an/9704001/}{funct-an/9704001}.
  \MR{2001g:05017}.

\bibitem{Knuth99a}
Donald~E. Knuth.
\newblock Convolution polynomials.
\newblock Preprint, Stanford Univ.,
  \href{http://www-cs-faculty.stanford.edu/~knuth/papers/cp.tex.gz}%
  {http://www-cs-faculty.stanford.edu/\~{}knuth/papers/cp.tex.gz}, February
  1996.

\bibitem{Rota95}
Joseph~P.S. Kung, editor.
\newblock {\em Gian-Carlo Rota on Combinatorics: Introductory Papers and
  Commentaries}, volume~1 of {\em Contemporary Mathematicians}.
\newblock Birkh\"auser Verlag, Boston, 1995.

\bibitem{Manin00a}
Yuri~I. Manin.
\newblock Classical computing, quantum computing, and {S}hor's factoring
  algorithm.
\newblock {\em Ast\'erisque}, (266):Exp.\ No.\ 862, 5, 375--404, 2000.
\newblock S\'eminaire Bourbaki, Vol. 1998/99. \MR{MR2001g:81040}.

\bibitem{Mattner99}
Lutz Mattner.
\newblock What are cumulants?
\newblock {\em Doc. Math.}, 4:601--622
  (\href{http://www.mathematik.uni--bielefeld.de/documenta/vol--04/18.ps.gz}{e%
lectronic}), 1999.
\newblock \MR{2001a:60017}.

\bibitem{Foundationiii}
Ronald Mullin and Gian-Carlo Rota.
\newblock On the foundation of combinatorial theory ({III}): Theory of binomial
  enumeration.
\newblock In B.Harris, editor, {\em Graph Theory and Its Applications}, pages
  167--213. Academic Press, Inc., New York, 1970.
\newblock Reprinted in~\cite[pp.~118--147]{Rota95}.

\bibitem{SimonReed75}
Michael Reed and Barry Simon.
\newblock {\em {Fourier} Analysis, Self-Adjointness}, volume~2 of {\em Methods
  of Modern Mathematical Physics}.
\newblock Academic Press, New York, {\noopsort{}}1975.

\bibitem{RomRota78}
S.~Roman and Gian-Carlo Rota.
\newblock The umbral calculus.
\newblock {\em Adv. in Math.}, 27:95--188, 1978.

\bibitem{Rota64a}
Gian-Carlo Rota.
\newblock The number of partitions of a set.
\newblock {\em Amer. Math. Monthly}, 71(5):498--504, May 1964.
\newblock Reprinted in~\cite[pp.~1--6]{Rota75} and~\cite[pp.~111--117]{Rota95}.

\bibitem{Rota75}
Gian-Carlo Rota.
\newblock {\em Finite Operator Calculus}.
\newblock Academic Press, Inc., New York, 1975.

\bibitem{KahOdlRota73}
Gian-Carlo Rota, David Kahaner, and Andrew Odlyzko.
\newblock Finite operator calculus.
\newblock {\em J. Math. Anal. Appl.}, 42(3):685--760, June 1973.
\newblock Reprinted in~\cite[pp.~7--82]{Rota75}.

\bibitem{Ryder96}
Lewis~H. Ryder.
\newblock {\em Quantum Field Theory}.
\newblock Cambridge University Press, Cambridge, 2nd edition,
  {\noopsort{1985}}1996.

\bibitem{ShabanovKlauder98}
Sergei~V. Shabanov and John~R. Klauder.
\newblock Path integral quantization and {R}iemannian-symplectic manifolds.
\newblock {\em Phys. Lett. B}, 435(3-4):343--349, 1998.
\newblock \MR{MR99i:81127}.
  \eprint{http://xxx.lanl.gov/abs/quant-ph/9805014}{quant-ph/9805014}.

\bibitem{Shubin87}
Mikhail~A. Shubin.
\newblock {\em Pseudodifferential Operators and Spectral Theory}.
\newblock Springer-Verlag, Berlin, 1987.

\bibitem{Simon79}
Barry Simon.
\newblock {\em Functional integration and quantum physics}.
\newblock Academic Press Inc. [Harcourt Brace Jovanovich Publishers], New York,
  1979.
\newblock \MR{MR84m:81066}.

\bibitem{MTaylor81}
Michael~E. Taylor.
\newblock {\em Pseudodifferential Operators}, volume~34 of {\em Princeton
  Mathematical Series}.
\newblock Princeton University Press, Princeton, New Jersey, 1981.

\end{thebibliography}
\end{document}